\documentclass[11pt]{amsart}
\usepackage{amssymb,mathrsfs,graphicx,enumerate,color}

\usepackage{colortbl}
\definecolor{black}{rgb}{0.0, 0.0, 0.0}
\definecolor{red}{rgb}{1.0, 0.5, 0.5}

\title[Stability versus Blow-up for Euler]{Blow-up solutions to 3D Euler are hydrodynamically unstable}

\author[Vasseur]{Alexis F. Vasseur}
\address[Alexis F. Vasseur]{\newline Department of Mathematics, \newline The University of Texas at Austin, Austin, TX 78712, USA}
\email{vasseur@math.utexas.edu}

\author[Vishik]{Misha Vishik}
\address[Misha Vishik]{\newline Department of Mathematics, \newline The University of Texas at Austin, Austin, TX 78712, USA}
\email{vishik@math.utexas.edu}


\thanks{\textbf{Acknowledgment.}  
The first author was partially supported by the NSF grant: DMS 1614918. }

\newcommand{\ve}{v_{\eps,\delta}}
\newcommand{\be}{\tilde{b}}
\newcommand{\C}{C_{\eta,\delta}}
\newcommand{\re}{r_{\eta,\delta}}
\newcommand{\Ve}{V_{\eps,\delta}}
\newcommand{\rre}{R_{\eps,\delta}}
\newcommand{\qe}{q_{\eps,\delta}}
\newcommand{\Ree}{\overline{R}_{\eps,\delta}}

\newcommand{\R}{\mathbb{R}}
\newcommand{\T}{\mathbb{T}}

\newcommand{\eps}{\varepsilon}
\newcommand{\dt}{\partial_t}
\newcommand{\Div}{\mathrm{div}\ }

\newcommand{\U}{{\Omega}}

\newcommand{\dis}{\displaystyle}

\newcommand{\dx}{\partial_x}

\newtheorem{theo}{Theorem}

\newtheorem{lemm}{Lemma}
\newtheorem{cor}{Corollary}
\newtheorem{prop}{Proposition}

\begin{document}
\bibliographystyle{plain}
\maketitle

\begin{abstract}
We   study the interaction  between the stability, and the  propagation of  regularity, for solutions to the incompressible 3D Euler equation. It is still unknown whether  a solution with smooth initial data can develop a singularity in finite time.  This article  explains why the  prediction of  such a blow-up, via direct numerical experiments, is so difficult. It is described how, in such a scenario,  the solution becomes unstable as time approaches the blow-up time.
\end{abstract}

\section{Introduction}

Consider  the incompressible Euler equation in a  domain $\U\subset \R^3$: 
\begin{equation} \label{1Euler}
\begin{array}{l}
\displaystyle{\dt u+ (u\cdot \nabla) u+\nabla P=0,\qquad x\in \U, \ \ t\in (0,T^*),}\\[0.3cm]
\displaystyle{\Div u=0, \qquad x\in \U, \ \ t\in (0,T^*),}
\end{array}
\end{equation}
endowed with  a  smooth initial value $u^0\in H^s(\U)$, for a $s>9/2$. The domain $\U$ can be $\R^3$, $\T^3$ or any smooth bounded domain $\U$ where we add the impermeability condition:
\begin{equation}\label{2Euler}
u\cdot \mathbf{n}=0, \qquad \mathrm{on} \ \ \partial \U, 
\end{equation}
where $\mathbf{n}$ is the normal of $\partial \U$.
 It is well known that there exists a  solution of this equation on (at least) a small timespan $(0,T^*)$ such that for every $T<T^*$, $u\in C^0(0,T;H^s(\U))\cap C^1(0,T; H^{s-1}(\U))$ (See for instance \cite{BKM}). From the assumption $s>9/2$, this implies that on this lifespan  $u$, $\nabla u$ and $\nabla\nabla u$ are $C^1([0,T]\times\U)$ for all $T<T^*$.  In two dimensions of space, due to the absence of vorticity stretching, the solution can always be extended as a global smooth solution for all time. Whether it is still the case in dimension 3, for all smooth initial values, remains one of the fundamental questions both for the Euler equation, and its viscous counterpart the Navier-Stokes equation.  This paper is dedicated to the study of the link between the linear stability of the solutions, and the propagation of their regularity.
Let  $T^*>0$ be the biggest time (possibly infinite) such that for every $T<T^*$ the solution of the Euler equation $u$ exists and lies in $C^0(0,T;H^s(\U))\cap C^1(0,T; H^{s-1}(\U))$. Let $1<p<\infty$.
For any $T<T^*$, we consider the semigroup generated by the linearization of the Euler equation about the solution $u$:

\begin{equation}\label{linear}
\begin{array}{ll}
&\dis{\dt v+ (u\cdot\nabla) v+(v\cdot\nabla) u+\nabla P'=0,\qquad x\in \U, \ \ t\in (0,T),}\\[0.2cm]
&\Div v=0, \qquad x\in \U, \ \ t\in (0,T),\\[0.2cm]
&v\cdot \mathbf{n}=0, \qquad \mathrm{on} \ \ \partial \U.
\end{array}
\end{equation}
The solution $v$ is uniquely defined for any initial value in $H^1(\U)$ (see Inoue and Miyakawa \cite{Inoue}). 
We denote $\gamma_p(T)$ the growth in  $L^p(\U)$  norm of the semi-group:
$$
\gamma_p(T)=\sup_{v^0\in H^1(\U), \|v^0\|_{ L^p(\U)}\leq1}\|v(T)\|_{L^p(\U)}.
$$
It is easy to show (see Lemma \ref{lemm2}) that the regularity on $[0,T]$ of the solution $u$ implies the boundedness of $\gamma_p(t)$ on $[0,T]$. Indeed, there exists a constant depending only on $\Omega$ and $p$ such that  
$$
\gamma_p(T)\leq e^{C_p \int_0^T\|\nabla u\|_{L^\infty(\U)}\,dt}.
$$
Therefore  regularity controls  linear stability. This paper is dedicated to the proof of  the other causality. 
We denote the vorticity $\omega= \mathrm{curl} \ u$.
Our main theorem is the following:
\begin{theo}\label{main_theo}
Let $\U\subset \R^3$ be either $\R^3$, $\T^3$, or a bounded smooth domain. Consider a smooth initial value $u^0\in H^s(\U)$, $s>9/2$, with $\omega^0=\mathrm{curl} \ u^0$, and denote $T^*$ the biggest time (possibly infinite) such that the solution  $u$ of the Euler equation (\ref{1Euler}) (\ref{2Euler}) exists and lies in $C^0(0,T;H^s(\U))\cap C^1(0,T; H^{s-1}(\U))$ for all $T<T^*$. Then, for every $1<p<\infty$,
and every $T< T^*$,
$$
\gamma^2_p(T)\geq  \frac{ \|\omega(T)\|_{L^\infty(\U)}}{\|\omega^0\|_{L^\infty(\U)}}.
$$
\end{theo}
Beale Kato and Majda showed in \cite{BKM} that the supremum norm of the vorticity controls the regularity of the Euler solution (see also Ferrari \cite{BKM2} in the bounded case).  
More precisely, they showed that, as long as 
$$
\int_0^T\|\omega(t)\|_{L^\infty(\U)}\,dt
$$
is bounded, there exists $\eps>0$ such that $u$ can be extended to a solution to the Euler equation (\ref{1Euler}) (\ref{2Euler}) on $[0,T+\eps]$ with $u\in C^0(0,T+\eps;H^s(\U))\cap C^1(0,T+\eps; H^{s-1}(\U))$.
Therefore, Theorem \ref{main_theo} implies the following result.

\begin{cor}\label{cor} Let $\U\subset \R^3$ be either $\R^3$, $\T^3$, or a bounded smooth domain. Consider a smooth initial value $u^0\in H^s(\U)$ for $s>9/2$ with $\omega^0=\mathrm{curl} \ u^0$, and assume that the corresponding solution $u$ to the Euler equation (\ref{1Euler}) (\ref{2Euler})  lies in $C^0(0,T;H^s(\U))\cap C^1(0,T; H^{s-1}(\U))$ for all $T<T^*$.
Assume that $$\sup_{0<T<T^*}\|u(T)\|_{H^s(\U)}=\infty. $$ Then, for any $1< p<\infty$:
$$
\int_0^{T^*} \gamma^2_p(t)\,dt=\infty,
$$
and especially  $$\lim\sup_{T\to T^*}\gamma_p(T)=\infty.$$
\end{cor}
This  is equivalent to the contrapositive which states that  stability controls the regularity. This result shows that if the solution $u$ blows up at $t=T^*$, then  small perturbations on the initial value induce  huge discrepancies on  the solution when time approaches $T^*$. Numerical experiments  involves unavoidable numerical inaccuracies. Therefore, due to the growing instabilities of the exact solution close to the blow-up time, we cannot expect any predictability of the numerical experiment about the blow-up. 
This  explains why, even with the current computational power,  it is so difficult to obtain numerical scenarios for a possible blow-up. The difficulty to predict finite time blow-ups is well documented (see  Hou and Li \cite{HouLi}, or Kerr \cite{Kerr} for instance). Note that the result of Corollary \ref{cor} covers the case of blow-ups at the boundary. Therefore it applies to the computation of Luo and Hou in  \cite{GuoHou}.
\vskip0.3cm
It is interesting to point out that this kind of situation is unusual. Consider compressible fluids for instance. Formation of singularities, known as shocks, are very easy to compute, and can be showed to be very stable. To illustrate the phenomenon, consider the simplified case of the one-dimensional scalar  Burgers equation:
$$
\dt u+\dx\left(\frac{u^2}{2}\right)=0.
$$
Shocks  happen in finite time, whenever the initial value is decreasing. 
The semigroup of the linearization of the Burgers equation about any solution $u$ is defined by the linear equation
$$
\dt v+\dx(uv)=0.
$$
In contrast  to the incompressible models, where the Energy norm $L^2$ comes naturally, the $L^1$ norm is more appropriate for the Burgers equation. Denote the growth of the $L^1$ norm of the semigroup as
$$
\rho_B(T)=\sup_{\|v_0\|_{ L^1(\R)}\leq 1}\|v(T)\|_{L^1(\R)}. 
$$
As long as $u$ is smooth enough, we have actually 
$$
\dt |v|+\dx(u|v|)=0.
$$
And so, integrating in $x$, we obtain that for all time up to the blowing-up time: $\rho_B(T)=1$. Hence the solution remains uniformly stable up to the formation of the shock. 
Note that in the case of the Burgers equation, any solutions is even  uniformly nonlinearly  stable in $L^1$ as proved by Kruzhkov \cite{K}.
\vskip0.3cm
The analogue  result for the Navier-Stokes equation is very easy (at least for $p>3/2$ and  $\U$ without boundary). Indeed, consider  $u^0\in L^2(\U)\cap W^{1,p}(\U)$, and $u$ a solution to the Navier-Stokes equation on $[0,T^*)\times\U$ with initial value $u^0$ and such that for  every $T<T^*$, $u\in L^\infty(0,T;W^{1,p}(\U))$. For every $i=1,2$ or 3, $v=\partial_{x_i}u$ is solution to the linearized Navier-Stokes equation:
\begin{equation*}
\begin{array}{ll}
&\dis{\dt v+ (u\cdot\nabla) v+(v\cdot\nabla) u+\nabla P'=\Delta v,\qquad x\in \U, \ \ t\in (0,T),}\\[0.2cm]
&\Div v=0, \qquad x\in \U, \ \ t\in (0,T).
\end{array}
\end{equation*}
Hence, the $L^p$ stability of the linearized Navier-Stokes equation  implies a bound in $L^\infty(0,T^*;L^p(\U))$ on $\nabla u$, which implies classical regularity on $u\in C^\infty((0,T^*+\eps)\times\U)$ for a $\eps>0$, as  long as $p>3/2$.
\vskip0.5cm
Let us now explain the general idea of the result. Let us say just for now that $\Omega=\R^3$. 
The vorticity $\omega=\mathrm{curl} \ u$ is solution to (see \cite{Constantin_Lecture} for instance):
$$
\dt \omega +(u\cdot\nabla) \omega-(\omega\cdot\nabla) u=0.
$$
But since 
\begin{equation}\label{omega}
((\nabla u)-(\nabla u)^T)\cdot\omega=0,
\end{equation} it is also solution to the following equation
$$
\dt \omega +u\cdot\nabla \omega-\nabla u\cdot\omega=0,
$$
which is dual to the linearized Euler equation (\ref{linear}). Therefore, for any backward solution $v$  to (\ref{linear}) with end-point value $v(T,\cdot)=v^T$,  we have
$$
\int_{\R^3} \omega(T,x)v^T(x)\,dx=\int_{\R^3} \omega^0(x)v(0,x)\,dx.
$$
If we denote  $\rho_b(T)$ the backward growth in  $L^1$  norm of the semi-group:
$$
\rho_b(T)=\sup_{\|v^T\|_{L^1(\R^3)}\leq 1}\|v(0)\|_{L^1(\R^3)},
$$
we get that 
$$
\|\omega(T)\|_{L^\infty(\R^3)} \leq \rho_b(T) \  \|\omega^0\|_{L^\infty(\R^3)}.
$$
Hence, from the Beale-Kato-Majda therorem \cite{BKM}, if 
$$
\int_0^{T^*}\rho_b(T)\,dT<\infty
$$
is finite, then the solution $u$ does not blow up at  time $T^*$. This provides a cheap proof that the backward stability of the linearized Euler equation in $L^1$ implies regularity. The difficulty is (1) to switch from backward stability to forward stability and (2) to extend the result to the  $L^p$ norms. Note that, even if the solution $u$ is reversible, the behavior at $t=0$, where $u$ is smooth,  is very different to the behavior at $t=T^*$, where $u$ is supposed to blow up. Therefore the  forward and backward linearization operators are very different.  Arnold showed  in \cite{Arnold} that, at least formally, the Euler equation has a natural symplectic structure. It can be therefore  expected that the forward stability growth is similar the backward one.  
\vskip0.3cm
To tackle this challenge, we use a WKB expansion method developed by the second author in   \cite{Vishik}  to study the essential spectral radius of small oscillations in an ideal incompressible fluid. The dynamic of the   solutions to (\ref{linear}) with high frequencies:
$$
v(t,\gamma_t(x_0))\approx b(t,x_0) e^{i S(t,x_0)/\eps} ,\qquad \mathrm{with} \ \ \nabla S(t,x_0)=\xi_t(x_0),
$$
for asymptotically  small $\eps$,  is described by the  bicharacteristic-amplitude equations, an ODE in time only on the unknown $(\gamma_t(x_0), \xi_t(x_0,\xi_0), b_t(x_0,\xi_0,b_0))$ where the dependence on $t$ is denoted as subscript, and $(x_0,\xi_0,b_0)$ are the initial values. The bicharacteristic-amplitude equations reads:
\begin{equation}\label{eq_algebric}
\left\{
\begin{array}{l}
\dot{\gamma_t}=u(t,\gamma_t), \\
\dot{\xi_t} =-(\partial_x u)^T \xi_t,\\
\dot{b_t}=-(\partial_x u)b_t+2\frac{\xi_t^T(\partial_x u)b_t}{|\xi|^2}\xi_t,
\end{array}
\right.
\end{equation}
where we denote $(\partial_x u)$ the matrix
$$
(\partial_x u)_{ij}=\partial_ju_i (t,\gamma_t).
$$

Note that  the incompressibility in the WKB expansion is expressed as
$$
b_t\cdot\xi_t=0,
$$
which is conserved by (\ref{eq_algebric}). The great advantage of studying the high frequencies solutions  via the WKB expansion is that several nonlocal properties of the linearized Euler equation become local. For instance, the (nonlocal) Leray projection involving the pressure, becomes local and corresponds to the second term on the right hand side of the equation on $b_t$.
The vorticity of the solution of the Euler equation along the trajectories $\omega_t(x_0) =\omega (t,\gamma_t(x_0))$ verifies
\begin{equation}\label{eq_omega} 
\dot{\omega_t}= (\partial_xu)\omega_t.
\end{equation}
We will also show that, as long as $\omega_0\cdot \xi_0=0$, the quantity $\omega_t\cdot b_t$ is conserved in time. This can be seen as the local conservation  (without  integration) of an helicity quantity involving the vorticity of the Euler equation, and the velocity of the linearized Euler equation.

\vskip0.3cm
The algebraic system of ODE (\ref{eq_algebric}) was used by Friedlander and the second author to define rigorously  the concept of  fluid Lyapunov exponent \cite{VishikDynamo}. They used it in several context to study the stability of steady states \cite{FriedVishik, FriedStrauss}.
\vskip0.3cm

In Section \ref{section2}, we perform a  thorough study of  the properties of this system of ODE (\ref{eq_algebric}). Let us denote
$$
\beta(T)=\sup_{|b_0|=1,|\xi_0|=1,  x_0\in \U, b_0\cdot\xi_0=0} |b_t(x_0,\xi_0,b_0)|.
$$
The quantity $\beta(T)$ is related to  the essential spectrum radius (see \cite{Vishik}).
We will show the following proposition.
\begin{prop}\label{prop1}
Let $\U$ be either $\R^3$, $\T^3$ or a bounded regular domain of $\R^3$. Assume that for a $s>9/2$, and every $T<T^*$, $u\in C^0(0,T;H^s(\U))\cap C^1(0,T;H^{s-1}(\U))$ and  is  solution on $[0,T^*)\times\U$ to the Euler equation (\ref{1Euler}) with, in the presence of boundary, the impermeability condition (\ref{2Euler}). Then, for every $T\in [0,T^*)$, we have:
\begin{equation}\label{vorticity-radius} 
\frac{1}{\sqrt{\|\omega^0\|_{L^\infty}}}\sup_{x\in \U}\sqrt{ |\omega(T,x)|} \leq \beta(T).
\end{equation}
\end{prop}
This shows that  the growth of $|\omega(t)|$ is at most as the square of the growth of $|b_t|$. 
\vskip0.5cm

In Section \ref{section3}, we will apply the WKB expansion to compare $\gamma_p(T)$ and $\beta(T)$. Namely, we will show the following.
\begin{prop}\label{prop2}
Let $\U$ be either $\R^3$, $\T^3$ or a bounded regular domain of $\R^3$. Assume that for a $s>9/2$, and every $T<T^*$, $u\in C^0(0,T;H^s(\U))\cap C^1(0,T;H^{s-1}(\U))$ and  is  solution on $[0,T^*)\times\U$ to the Euler equation (\ref{1Euler}) with, in the presence of boundary, the impermeability condition (\ref{2Euler}). 
Then, for any $T<T^*$ and any $1<p< \infty$, we have
$$
\beta(T)\leq \gamma_p(T).
$$
\end{prop}

\vskip0.5cm
Theorem \ref{main_theo}, is a direct application of Proposition \ref{prop1} and Proposition \ref{prop2}. Indeed, Assume that for a $s>9/2$, and every $T<T^*$, $u\in C^0(0,T;H^s(\U))\cap C^1(0,T;H^{s-1}(\U))$ and  is  solution on $[0,T^*)\times\U$ to the Euler equation (\ref{1Euler}) with, in the presence of boundary, the impermeability condition (\ref{2Euler}). Then applying both propositions, we obtain that for every $1<p<\infty$:
$$
\sup_{x\in \U} |\omega(T,x)| \leq (\beta(T))^2 \|\omega^0\|_{L^\infty(\U)}\leq (\gamma_p(T))^2 \|\omega^0\|_{L^\infty(\U)},
$$
which is the statement of Theorem \ref{main_theo}.

\section{Vorticity and  essential spectrum radius}\label{section2}

This section is dedicated to the proof of Proposition \ref{prop1}. 
For $t$ fixed, consider $\gamma_t:\U\to\U$. When the function is invertible, we denote $\gamma_t^{-1}$ his inverse. That is, the function such that for every $x_0\in \Omega$ $\gamma_t(\gamma_t^{-1}(x_0))=x_0$. 
We begin with a lemma. 
\begin{lemm}\label{lemm1}
Let $\U$ be either $\R^3$, $\T^3$ or a bounded regular domain of $\R^3$.  for a $s>9/2$, and every $T<T^*$, $u\in C^0(0,T;H^s(\U))\cap C^1(0,T;H^{s-1}(\U))$ and  is  solution on $[0,T^*)\times\U$ to the Euler equation (\ref{1Euler}) with, in the presence of boundary, the impermeability condition (\ref{2Euler}).  Let $T<T^*$. Then, for any  $x_0\in \U$, and  $\xi_0, b_0, \be_0\in \R^3$,  there exist unique solutions  $(\gamma_t,\xi_t, b_t)$ and $(\gamma_t,\xi_T, \be_t)$  to (\ref{eq_algebric}) with initial values $(x_0,\xi_0,b_0)$, and $(x_0,\xi_0,\be_0)$.  For every $t\in [0,T]$, the function $\gamma_t$ is invertible. The functions $\gamma, \gamma^{-1}, \xi, b, \be$ are Lipschitz with respect to time and $C^2$ with respect to the initial values.  For all time $t\in [0,T]$, $\gamma_t\in \U$. In addition, if  $b_0, \be_0, \xi_0$ are linearly independent,  $(b_0\cdot \xi_0)=(\be_0\cdot \xi_0)=0$,  and $\omega_t$ is solution to (\ref{eq_omega}), then  we have:
\begin{eqnarray}
\label{cons_omega}
&&\omega_T\cdot\xi_T=\omega_0\cdot \xi_0,\\
\label{cons_b}
&& b_T\cdot\xi_T=\be_T\cdot \xi_T=0,\\
\label{cons_B}
&&(b_T\times \be_T)\cdot \xi_T=(b_0\times \be_0)\cdot\xi_0.
\end{eqnarray}
\end{lemm}
Note that $\gamma$ does not depend on $\xi$ nor on $b$, while $\xi$ does not depend on $b$. We can then construct $b$ and $\be$ for the same functions $\gamma$, and $\xi$.
\begin{proof}
First, note that by Sobolev embedding,  $u$, $\nabla u$ and $\nabla\nabla u$ lie  in $C^1([0,T]\times \U)$. Hence Cauchy-Lipschitz provides a unique solution to (\ref{eq_algebric}) which is Lipschitz both in time and with respect to the initial conditions. The map $\gamma_t^{-1}$ can be defined as $\bar\gamma_t$ using the flow $\bar{u}(s,x)=-u(t-s,-x)$. Hence $\gamma_t^{-1}$ is also Lipschitz in both time and $x$.
Differentiating (\ref{eq_algebric}) with respect to the initial values, shows that those quantities are $C^2$ with respect to the initial conditions. Moreover, for  $x_0\in \U$ open, for all $t\in [0,T]$ , $\gamma_t\in \U$. Indeed, if for a $t_0$, $\gamma_{t_0}\in \partial\U$, then, from the permeability condition, $\gamma_t\in \partial\U$ for all $t\in [0,T]$, which is not possible since $x_0\in \U$. 
\vskip0.5cm
The equation (\ref{eq_omega}) and the second equation of (\ref{eq_algebric}) are dual, so (\ref{cons_omega}) is straightforward. 
Similarly, we obtain directly from the second equation and the third equation of (\ref{eq_algebric}) that $ b\cdot\xi$ is conserved by the flow, and so equal to 0 thanks to the initial values. Note that this corresponds to the incompressibility of the Euler linearization for high  frequencies.
\vskip0.5cm
For three vectors $B_1,B_2,B_3\in \R^3$ let us denote   $(B_1,B_2,B_3)$ the matrix $(B_{ij})$ where $B_{ij}$ is the $i$th component of the vector $B_j$. We denote the matrix
$$
\Psi_t=(b_t,\be_t,\xi_t).
$$
We have 
$$
\dot{\Psi_t}=-(\partial_x u)\Psi_t+(\alpha_1(t) \xi_t, \alpha_2(t)\xi_t, \xi_t\times\omega),
$$
where we used
\begin{eqnarray*}
&&\alpha_1(t) =2\frac{\xi_t^T(\partial_xu) b_t}{|\xi_t|^2},\\
&&\alpha_2(t) =2\frac{\xi_t^T(\partial_xu) \be_t}{|\xi_t|^2},\\
&& -(\partial_x u)^T\xi_t=-(\partial_x u)\xi_t+[(\partial_x u)-(\partial_x u)^T]\xi_t=-(\partial_x u)\xi_t+\omega\times\xi_t.
\end{eqnarray*}
We want to compute the evolution of 
$$
(b\times \be)\cdot\xi=-\mathrm{Det}(\Psi),
$$
where we dropped the subscript expressing  the dependence in $t$ to simplify the notations.
As long as $\Psi$ is invertible, we  have 
\begin{eqnarray*}
&&\qquad\qquad\frac{d}{dt}{\mathrm{Det}(\Psi)}=\mathrm{Tr} (\Psi^{-1}\dot{\Psi}) \mathrm{Det}(\Psi)\\
&&=\left[-\mathrm{Tr} (\Psi^{-1}(\partial_x u)\Psi) + \mathrm{Tr}  (\alpha_1 \Psi^{-1}\xi, \alpha_2\Psi^{-1}\xi,\Psi^{-1} [\xi\times\omega]))\right]\mathrm{Det}(\Psi).
\end{eqnarray*}
Since $\mathrm{Tr} (AB)=\mathrm{Tr}(BA)$,
$$
\mathrm{Tr} (\Psi^{-1}(\partial_x u)\Psi)=\mathrm{Tr} ((\partial_x u)\Psi\Psi^{-1})=\mathrm{Tr}(\partial_xu)=\mathrm{div} \ u=0.
$$
Note that by definition, $\Psi e_1=b$, $\Psi e_2=\be$, and $\Psi e_3=\xi$. Hence
$$
 \Psi^{-1}b=e_1, \qquad \Psi^{-1}\be=e_2, \qquad \mathrm{and} \ \ \Psi^{-1}\xi=e_3.
$$
Since $\xi\times\omega$, $b$ and $\be$ are all orthogonal to $\xi$, and $b,\be, \xi$ are linearly independent,  there exists scalars $\beta_1$ and $\beta_2$ depending only on time such that
$$
\xi\times\omega=\beta_1 b+\beta_2 \be.
$$
We denote $E_{ij}$ the matrix such that $(E_{ij})_{kl}=\delta_{kl}$. Note that $\mathrm{Tr} E_{ij}=\delta_{ij}$.
We have 
\begin{eqnarray*}
&&\mathrm{Tr}  (\alpha_1 \Psi^{-1}\xi, \alpha_2\Psi^{-1}\xi,\Psi^{-1} [\xi\times\omega]))=\mathrm{Tr}(\alpha_1 e_3, \alpha_2 e_3, \beta_1e_1+\beta_2 e_2 )\\
&&\qquad =\alpha_1 \mathrm{Tr} E_{31}+\alpha_2\mathrm{Tr} E_{32}+\beta_1\mathrm{Tr} E_{13} +\beta_2 \mathrm{Tr} E_{23}=0.
\end{eqnarray*}
This shows that $\mathrm{det} (b,\be,\xi)$ is constant in time. Hence the matrix $\Psi$ remains invertible, and the computation makes sense for all $t\in [0,T]$. This provides
(\ref{cons_B}). Note that this result is due to the incompressibility of the flow generated by $u$, and represents the conservation of volume by the flow. 
\end{proof}
We can now show  Proposition \ref{prop1}. 
\begin{proof}
Fix $(T,x_T)\in (0,T^*)\times\U$. Define (backward) the solution $\gamma_t$ of the first equation of (\ref{eq_algebric}) such that $\gamma_T(x_0)=x_T$. As in Lemma \ref{lemm1}, $\gamma_t$ is uniquely defined on $[0,T]$ and for all $t$, $\gamma_t\in \U$. Consider $\xi_T$ the unit vector in the direction of $\omega(T,x_T)$ such that $\xi_T\cdot\omega(T,x_T)=|\omega(T,x_T)|$. From the first equality of Lemma \ref{lemm1},
\begin{equation}\label{1}
|\omega(T,x_T)|=\xi_0\cdot \omega^0(x_0)\leq  |\xi_0|\sup_\U|\omega^0|.
\end{equation}
Consider unit vectors $b_0, \be_0$ such that $b^0, \be_0, \xi_0$ are orthogonal to each others, and $(b_0\times \be_0)\cdot\xi_0=|\xi_0|$. Consider $b, \be$ the solutions of the third equation  of (\ref{eq_algebric}) corresponding to the same $\xi_t$, and $\xi_t$ constructed above. From the last equality of Lemma \ref{lemm1}, we have 
\begin{equation}\label{eq_12}
 |\xi_0|=(b_0\times \be_0)\cdot\xi_0=(b_T\times \be_T)\cdot\xi_T\leq |b_T||\be_T|\leq \left(\sup_{|b_0|=1,\xi_0\in\R,  x_0\in \U, b_0\cdot\xi_0=0} |b_T|\right)^2.
\end{equation}
Note that if $(\gamma, \xi,b)$ is solution to (\ref{eq_algebric}), then $(\gamma, \xi/|\xi_0|, b)$ is also solution to (\ref{eq_algebric}). Therefore,
$$
\sup_{|b_0|=1,\xi_0\in\R,  x_0\in \U, b_0\cdot\xi_0=0} |b_T(x_0,\xi_0,b_0)|=\sup_{|b_0|=1,|\xi_0|=1,  x_0\in \U, b_0\cdot\xi_0=0} |b_T(x_0,\xi_0,b_0)|,
$$
and together with (\ref{eq_12}) and (\ref{1}), this gives (\ref{vorticity-radius}).
\end{proof}

\section{Linear stability  and essential spectrum radius}\label{section3}

This section is dedicated to the proof of Proposition \ref{prop2}. We begin with some  some easy results on the linearized Euler equation with source term:
\begin{equation}\label{eq_lin}
\begin{array}{l}
\dis{\dt v+ (u\cdot\nabla) v+(v\cdot\nabla) u+\nabla P'=f,\qquad x\in \U, \ \ t\in (0,T),}\\[0.3cm]
\dis{\Div v=0, \qquad x\in \U, \ \ t\in (0,T),}\\[0.3cm]
\dis{v\cdot \mathbf{n}=0, \qquad \mathrm{on} \ \ \partial \U,}\\[0.3cm]
\dis{v(0,\cdot)=v_0\qquad x\in \U.}
\end{array}
\end{equation}
\begin{lemm}\label{lemm2}
Let $\U$ be either $\R^3$, $\T^3$ or a bounded regular domain of $\R^3$. Assume that for a $s>5/2$, and every $T<T^*$, $u\in C^0(0,T;H^s(\U))\cap C^1(0,T;H^{s-1}(\U))$ and  is  solution on $[0,T^*)\times\U$ to the Euler equation (\ref{1Euler}) with, in the presence of boundary, the impermeability condition (\ref{2Euler}). 
 Let $1<p<\infty$, then there exists a constant $C_p>0$ such that for every $T<T^*$, and   
$v^0\in H^1(\U)\times L^p(\U)$  and $f\in L^2(0,T; H^1(\U))\cap L^1(0,T; L^p(\U))$, the following is true. There exists a unique solution $v\in C^0(0,T; H^1(\U))$ to the Linearized equation (\ref{eq_lin}) 
 with the impermeability condition $v\cdot \mathbf{n}=0$ at the boundary $\partial\U$. Moreover :
\begin{eqnarray*}
&&\| \nabla P'(t)\|_{L^p(\U)}\leq C_p (\|u(t)\|_{H^s(\U)} \|v(t)\|_{L^p(\U)}+\|f(t)\|_{L^p(\U)}), \qquad \mathrm{for \ \ all \ \  } t\in [0,T],\\
&&\|v(T)\|_{L^p(\U)}\leq e^{C_p \int_0^T\|u(t)\|_{H^s(\U)}\,dt}\left( \|v^0\|_{L^p(\U)}+C_p\|f\|_{L^1(0,T;L^p(\U))}\right).
\end{eqnarray*}
\end{lemm}
Especially, it shows that the regularity of $u$ implies the linear stability of the solution in any $L^p$, $1<p<\infty$. Note that, in the case of boundary, the regularity on the initial value $v^0\in H^1$ provides enough regularity on $v(t)$ to make sense of the impermeability condition on $v$.   
\begin{proof}
Since $v^0\in H^1(\U)$, and $u\in C^1((0,T)\times \U)$, Inoue and Miyakawa \cite{Inoue} insures the existence and uniqueness of the solution in $C^0(0,T; H^1(\U))$ verifying the impermeability condition at the boundary.   Taking the divergence of the first equation (\ref{linear}) gives:
\begin{eqnarray*}
&&-\Delta P'=2\Div[ (\partial_x u)v-f],\qquad \mathrm{in}  \ \ \U.
\end{eqnarray*}
In the presence of a boundary $\partial \U$, we can decompose any vector field $V$ on $\partial \U$ as $V= V_\tau+(V\cdot\mathbf{n})\mathbf{n}$. The impermeability conditions on both $u$ and $v$ implies that on the boundary $u=u_\tau$ and $v=v_\tau$. Hence, for any $x\in\partial\U$
 $$
 [(\partial_x u)v]\cdot \mathbf{n}=[(v_\tau\cdot \nabla_\tau) u]\cdot\mathbf{n}= (v_\tau\cdot \nabla_\tau)(u \cdot \mathbf{n})=0.
 $$
 Similarly, 
 $$
  [(\partial_x v)u]\cdot \mathbf{n}=0, \qquad \mathrm{on} \ \ \partial\U.
 $$
 Hence, from the equation, $\partial_n P'=f\cdot \mathbf{n}$ on $\partial \U$, which can be rewritten:
 $$
- \partial_n P'= [2(\partial_x u)v -f]\cdot \mathbf{n}  \qquad \mathrm{on} \ \ \partial\U.
 $$
Therefore, thanks to  classical elliptic regularity (see  Krylov \cite{Krylov} fo instance), we have the existence of $C_p>0$ such that:
\begin{eqnarray*}
&&\| \nabla P'(t)\|_{L^p(\U)}\leq C_p (\|\partial_x u(t)\|_{L^\infty(\U)} \|v(t)\|_{L^p(\U)}+\|f(t)\|_{L^p(\U)})\\
&&\qquad\qquad \leq C_p (\| u(t)\|_{H^s(\U)} \|v(t)\|_{L^p(\U)}+\|f(t)\|_{L^p(\U)}),
\end{eqnarray*}
thanks to the Sobolev embedding since $s>5/2$. 
Then multiplying the linearized Euler equation by $p|v|^{p-2}v$, integrating in $x$, and dividing by $\|v\|^{p-1}_{L^p}$ gives
\begin{eqnarray*}
&&\frac{d}{dt}\|v(t)\|_{L^p(\U)}\leq p\|\partial_x u(t)\|_{L^\infty(\U)} \|v(t)\|_{L^p(\U)}+p \|P'(t)\|_{L^p(\U)}+p\|f(t)\|_{L^p(\U)}\\
&&\qquad\qquad \leq  p(C_p+1)( \|u(t)\|_{H^s(\U)} \|v(t)\|_{L^p(\U)}+\|f(t)\|_{L^p(\U)}).
\end{eqnarray*}
Using the Gronwall argument provides the second estimates (for a bigger constant $C$ that we relabeled $C_p$). 
\end{proof}
We now show Proposition \ref{prop2}.

\begin{proof}
From Proposition \ref{prop1}, $\gamma_t$, $\gamma_t^{-1}$, $b_t$ and $\xi_t$ are Lipschitz functions  in time and $C^2$  with respect to the initial values $(x_0, \xi_0, b_0)$. The regularity control  depends only on $u^0$, $\U$, and  $T$. Let $T\in [0,T^*)$.
For any $\eta>0$, 
 Consider $x_0\in \U$, $\xi_0$ and $b_0$ with $|b_0|=|\xi_0|=1$ such that 
$$
\beta(T) \leq (1+\eta) |b_T(x_0,\xi_0, b_0)|.
$$
There exists $\delta>0$ such that the ball centered at $\gamma_T(x_0)$ with radius $\delta$ is strictly in $\U$ and 
$$
\inf_{x\in B_\delta (\gamma_T(x_0))} |b_T(g_T^{-1}(x),\xi_0, b_0)| \geq (1-\eta) |b_T(x_0, \xi_0, b_0)|. 
$$
Let $\overline{\phi}$ be a  a regular nonnegative function compactly supported in $B_1$ with $L^p$ norm equal to 1, and  
$$
\phi_T(x)=\frac{1}{\delta^{3/p}}\overline{\phi}\left(\frac{x-x(T)}{\delta}\right).
$$
We have $\|\phi_T\|_{L^p}=1$.
For any $\eps>0$, with a small abuse of notation, we denote the following functions of $(t,x)$ defined on $[0,T]\times\U$:
\begin{eqnarray*}
&&\xi(t,x)=\xi_t(g^{-1}_t(x), \xi_0),\\
&& b(t,x)=b_t(g^{-1}_t(x), \xi_0,b_0),\\
&&S(t,x)=(g^{-1}_t(x))\cdot \xi_0,\\
&& \phi(t,x)=\phi_T(g_t^{-1}(x)).
\end{eqnarray*}
We  consider the potential defined on $[0,T]\times\U$ as
$$
A_{\eps,\delta}=\eps \frac{\xi\times b}{|\xi|^2} \phi e^{i\frac{S}{\eps}},
$$
and 
$$
\ve=\mathrm{curl} \ A_{\eps,\delta}.
$$
obviously 
$$
\Div \ve=0.
$$
The function $\phi(T,\cdot)$ is compactly supported in $\U$, therefore, thanks to Lemma \ref{lemm1}, $\ve(t,\cdot)$ is compactly supported in $\U$ for all $t\in[0,T]$.
Hence 
$$
 \ve\cdot\nu=0,\qquad \mathrm{on} \ \partial\U.
$$

The function $S$ and $\phi$ verify the transport equation:
\begin{eqnarray*}
&& \dt \phi+ (u\cdot\nabla) \phi=0,\qquad t,x\in (0,T)\times\U,\\
&& \phi(T,x)=\Phi_T(x), \qquad x\in \U,
\end{eqnarray*}
\begin{eqnarray*}
&& \dt S+ (u\cdot\nabla) S=0,\qquad t,x\in (0,T)\times\U,\\
&& S0,x)=x\cdot\xi_0, \qquad x\in \U.
\end{eqnarray*}
Then 
$
i\phi  e^{i S/\eps}
$
is also transported by the flow, and verifies the same equation.
We denote 
$$
\Ve=i\phi b e^{i S/\eps}.
$$
\vskip0.3cm
The function $\nabla S$ verifies the equation:
$$
 \dt \nabla S+ (u\cdot\nabla)\nabla S=-(\partial_x u)^T\nabla S,\qquad t,x\in (0,T)\times\U,
$$
with initial value $\xi_0$.
Hence $\nabla S(t, g_t(x_0))$ is solution to the second equation of (\ref{eq_algebric}), and so
 \begin{eqnarray*}
 &&\nabla S(t,x)=\xi(t,g^{-1}_t(x),\xi_0)=\xi(t,x).
 \end{eqnarray*}
 \vskip0.5cm
 We fix $\eta$ and $\delta$.
 All the function being $C^1$ in time and $C^2$ in $x$  on $[0,T]\times \U$, we want to track the dependence on $\eps$  only.
 Note that the dependence  on $\eps$ comes only form the function $\eps e^{iS(t,x)/\eps}$. In the following we denote $\C$ any constant which does not depend on $\eps$. 
Since $\xi_0\cdot b_0=0$, from Lemma \ref{lemm1}, $\xi\cdot b=0$ on $[0,T]\times\U$, so
\begin{eqnarray*}
&&-\frac{i}{\eps} \nabla S\times A_{\eps,\delta}=-i \frac{\xi\times(\xi\times b)}{|\xi|^2} \phi e^{i S/\eps}\\
&&\qquad\qquad = i \frac{b(\xi\cdot\xi)-\xi(\xi\cdot b)}{|\xi|^2} \phi e^{i S/\eps}\\
&&\qquad\qquad = i \phi b  e^{i S/\eps}=\Ve.
\end{eqnarray*}
therefore we have $\ve=\Ve+\eps e^{iS/\eps}\re$ where
$$
\re=\mathrm{curl} \left( \frac{\xi\times b}{|\xi|^2} \phi\right)
$$
does not depend on $\eps$ and so, for all $t\in [0,T]$:
$$
\|\re(t)\|_{L^p}\leq \C.
$$
 Since $e^{iS/\eps}$ is transported by the flow, we get
 \begin{equation}\label{*}
 \|(\dt  +(u\cdot\nabla))(\eps e^{iS/\eps}\re)\|_{L^p}\leq \C\eps. 
\end{equation}
Moreover $\Ve$ verifies
\begin{eqnarray*}
&&\dt \Ve+(u\cdot\nabla)\Ve=\phi e^{iS/\eps} ( \dt b+ (u\cdot\nabla) b)\\
&&\qquad =-\phi e^{iS/\eps} \left((\partial_x u)b-2\xi\frac{\xi^T(\partial _xu)b}{|\xi|^2}\right)\\
&&=-(\Ve\cdot\nabla)u-2\left(\frac{\xi^T(\partial_x u)\Ve}{|\xi|^2}\right)\xi.
\end{eqnarray*}
 Denote 
 $$
 \qe=-2i\eps \frac{\xi^T(\partial_x u)\Ve}{|\xi|^2}.
 $$
 We have 
 $$
 -\nabla \qe=2\left(\frac{\xi^T(\partial_x u)\Ve}{|\xi|^2}\right)\xi +\rre,
 $$
 with 
 $$
 \|\rre\|_{L^p}\leq \C\eps, \qquad \mathrm{for \ all \ } t\in [0,T].
 $$
Together with (\ref{*}), this gives
 $$
 \dt \ve+(u\cdot\nabla)\ve+(\ve\cdot \nabla)u+\nabla q_{\eps,\delta}=\Ree,
 $$
 with 
 $$
 \|\Ree(t)\|_{L^p}\leq \C\eps, \qquad \mathrm{for \ all \ } t\in [0,T].
 $$
 Consider $v$ the solution to (\ref{linear}) with the same initial value as $\ve$. Applying Lemma \ref{lemm2} to $v-\ve$, we find that
$$
\|v(T)-\ve(T)\|_{L^p(\U)}\leq C\|\Ree\|_{L^1(0,T;L^p(\U))}\leq \C \eps.
$$ 
Hence
 \begin{eqnarray*}
&& \gamma_p(T)\geq \|v(T)\|_{L^p(\U)} \geq \|\ve(T)\|_{L^p(\U)} -\|v(T)-\ve(T))\|_{L^p(\U)} \\
 \qquad &&\qquad \geq \|\Ve(T)\|_{L^p(\U)} - \C \eps\\
  \qquad &&\qquad \geq \|\phi_T b_T)\|_{L^p(\U)} - \C \eps\\
 \qquad &&\qquad\geq (1-\eta) |b(T,x_0,\xi_0)| - \C \eps\\
 \qquad &&\geq \frac{1-\eta}{1+\eta} \beta (T) - \C \eps.
 \end{eqnarray*} 
 Taking the limits $\eps$ goe to zero,  then $\eta$ goes to zero gives the result.
 \end{proof}

\bibliography{biblio.bib}

\end{document}